\documentclass[a4paper,11pt,oneside]{amsart}
\usepackage{geometry}
\usepackage{graphicx, amsfonts, amssymb}
\usepackage{mathrsfs, amsmath, amsthm}
\usepackage{verbatim}
\usepackage{relsize}
\usepackage{enumerate}
\usepackage{hyperref}
\usepackage{version}
\usepackage{blindtext}
\usepackage{mathtools}
\usepackage{color}

\def\bbC{{\mathbb C}}
\def\bbD{{\mathbb D}}

\def\bbN{{\mathbb N}}

\def\bbR{{\mathbb R}}

\def\bbT{{\mathbb T}}

\def\bbZ{{\mathbb Z}}

\def\cF{{\mathcal F}}

\def\cH{{\mathcal H}}

\def\cP{{\mathcal P}}

\def\cS{{\mathcal S}}

\newcommand{\ra}{\rangle}
\newcommand{\la}{\langle}

\def\bpm{\begin{pmatrix}}
\def\epm{\end{pmatrix}}

\def\epf{ $\Box$\medskip

}

\def\ms{\medskip}

\def\vp{\varphi}

\def\ov{\overline}

\def\tbi{\item[{\tiny $\bullet$}]}

\def\supp{\operatorname{supp}}

\DeclareMathOperator{\Hol}{Hol}

\DeclareMathOperator{\BMOA}{BMOA}
\DeclareMathOperator{\VMOA}{VMOA}
\DeclareMathOperator{\BMO}{BMO}
\DeclareMathOperator{\LMO}{LMO}

\hypersetup{
	colorlinks=true,
	linkcolor=blue,
	filecolor=magenta,      
	urlcolor=cyan,
	pdftitle={Meromorphic Optimal Domain of Integral Operators},
	pdfpagemode=FullScreen,
}

\newtheorem{theorem}{Theorem}[section]
\newtheorem{Lemm}[theorem]{Lemma}
\newtheorem{prop}[theorem]{Proposition}

\newtheorem{openproblem}{Question}

\theoremstyle{definition}
\newtheorem{definition}[theorem]{Definition}

\theoremstyle{remark}

\numberwithin{equation}{section}

\begin{document}

\title{On Toeplitz operators on $H^1(\bbC^+)$}
\author{Carlo Bellavita}
\address{Departament de Matem\'atiques i Inform\'atica, Universitat de Barcelona,
Gran Via 585, 08007 Barcelona, Spain}
\email{carlobellavita@ub.edu}

\author{Marco M. Peloso}
\address{Dipartimento di Matematica F. Enriques, Dipartimento di Eccellenza MUR 2023-2027, Universit\`a degli Studi di Milano, Via C. Saldini 50, I-20133 Milano}
\email{marco.peloso@unimi.it}

\thanks{The authors are members of Gruppo Nazionale per l’Analisi Matematica, la Probabilit\`a e le loro Applicazioni (GNAMPA) of Istituto Nazionale di Alta Matematica (INdAM).
The first author was partially supported by PID2021-123405NB-I00 by
the Ministerio de Ciencia e Innovaci\'on. He wishes to thank for the support provided.
The second author was partially supported by the 2022 INdAM–GNAMPA project Holomorphic Functions in One and Several Complex Variables
(CUP E55F22000270001). He wishes to thank GNAMPA for the support provided.}

\date{\today}



\begin{abstract}
In this paper we consider Toeplitz operators with anti-analytic
symbols on $H^1(\bbC^+)$.  It is well known
that there are no bounded
Toeplitz operators
$T_{\ov\Theta}\colon H^1(\bbC^+)\to H^1(\bbC^+)$, where $\Theta\in H^\infty(\bbC^+)$. We consider the
subspace $H^1_\Theta=\{f\in H^1(\bbC^+)\colon \int_\bbR f\ov\Theta =0\}$
and show that it is natural to study the boundedness of
$T_{\ov\Theta}\colon H^1_\Theta \to H^1(\bbC^+)$.  We provide several
different conditions equivalent to such boundedness. We
prove that when $\Theta=e^{i\tau(\cdot)}$, with $\tau>0$,
$T_{\ov\Theta}\colon H^1_\Theta \to H^1(\bbC^+)$ is bounded. Finally,
we discuss a number of related open questions.
\end{abstract}

\maketitle

\section*{Introduction}

Let $\bbD$ be the unit disk and $\bbC^+$ the upper half-plane in
$\bbC$. For $p\in (0,\infty]$, let $H^p(\bbD)$ and $H^p(\bbC^+)$ denote the 
Hardy spaces on $\bbD$ and $\bbC^+$, resp.
We write $H^p$ to denote either the space on $\bbD$ or $\bbC^+$.
Functions in $H^p$ admits non-tangential boundary values on $\bbT$ or on
$\bbR$, whichever is appropriate. Given a function $f\in H^p$, as customary, we denote by $f$
also its boundary-value function.  The mappings $H^p(\bbD)\ni f\mapsto f\in L^p(\bbT)$,
and $H^p(\bbC^+)\ni f\mapsto f\in L^p(\bbR)$ are isometries, where the
measure on $\bbT$ is the normalized Lebesgue measure.  Thus, functions
in $H^p$ can be identified with their 
 boundary value functions.  
We also consider the spaces of anti-analytic functions $\ov{H^p}$,
both in $\bbD$ and $\bbC^+$.  Notice that $H^p(\bbC^+)\cap
\ov{H^p}(\bbC^+) =\{0\}$, while $H^p(\bbD)\cap
\ov{H^p}(\bbD) =\operatorname{span}\{1\}$.

When $p=2$, $H^2$  is a
reproducing kernel Hilbert space, and also a closed subspace of
$L^2(\bbT)$ and $L^2(\bbR)$, resp.  We denote by $P_+$ the orthogonal
projection from $L^2$ onto $H^2$, and by $P_-$ the orthogonal
projection from $L^2$ onto $\ov{H^2}$, in both the cases of $\bbD$ and
$\bbC^+$.

The {\em Toeplitz operator}
with symbol $\vp$ is defined as
$$
T_\vp \colon H^2\ni f \mapsto P_+(\vp f)\in   H^2 ,
$$
whenever it is well defined, while the {\em Hankel operator} with
symbol $\vp$ is defined as
$$
H_\vp \colon H^2\ni f \mapsto P_-(\vp f)\in  \ov{ H^2 }.
$$ 
It is well known that $T_\vp\colon H^2\to H^2$ is bounded if and only if
$\vp\in L^\infty$, while $H_\vp \colon H^2\to \ov{H^2}$ is bounded if and only
if $\vp\in \BMO$, see e.g.~\cite[Part B]{Nikolski} for these results, and also 
Section~\ref{Sec1} for the definition of $\BMO$. 
The above mentioned results extend to the case of $H^p$, for
$p\in(1,\infty)$,   see \cite{Janson1984}. 
In the case $p=1$ on the unit disk, 
S.\ Janson, J.\ Peetre and S.\ Semmes in
\cite{Janson1984} showed that 
the Toeplitz operator $T_{ \vp}\colon H^1(\bbD)\to H^1(\bbD)$
is bounded if and only if $\vp$ belongs to $L^\infty(\bbT) 
\cap \LMO(\bbT)$, where $\LMO(\bbT)$ is 
a $\log$-$\BMO$ type space, see Section~\ref{2:sec}.  In the same
article  Janson,  Peetre and  Semmes pointed out that in the case of
$\bbC^+$ there are no bounded Toeplitz operator
$T_{\ov\Theta}\colon H^1(\bbC^+)\to H^1(\bbC^+)$,  when $\Theta \in H^\infty$.

In this paper we consider the  problem  of the boundedness of Toeplitz operators on 
$H^1(\bbC^+)$. We show that the correct question is to consider the
Toeplitz operators defined only on a (closed) subspace of
$H^1(\bbC^+)$, defined in the terms of the symbol of the  Toeplitz
operator. 
Then, we restrict our attention to anti-analytic
symbols $\ov\Theta$, where $\Theta\in H^\infty(\bbC^+)$ and
we denote by $H^1_\Theta$ the closed subspace of $H^1(\bbC^+)$, given
by $H^1_\Theta=\{f\in H^1(\bbC^+)\colon \int f \ov\Theta=0\}$, to
which we restrict the Toeplitz operator $T_{\ov\Theta}$.  We provide
some conditions that are equivalent to the boundedness of
$T_{\ov\Theta}\colon H^1_\Theta \to H^1(\bbC^+)$, Theorem
\ref{equiv-cond:thm}.  The final main result is Theorem
\ref{main-result} in which we prove that $T_{\ov\Theta}\colon
H^1_\Theta \to H^1(\bbC^+)$ is bounded when $\Theta=e^{i\tau(\cdot)}$,
with $\tau>0$. 
\ms

\section{Basic facts and statement of the main results}\label{Sec1}
For $0< p \leq \infty$, the Hardy space $H^p(\bbC^+)$ of the
upper half-plane  $\bbC^+$ consists of  the analytic
functions $f$ on $\bbC^+$ such that, when $0<p<\infty$, 
$$
\|f\|_{H^p}^p = \sup_{y > 0} \int_{\bbR} |f(x+iy)|^p dx < \infty
$$
and, when $p=\infty$, 
$$
\|f\|_{H^\infty} = \sup_{z \in \bbC^+} |f(z)|<\infty.
$$
If $1\le p<\infty$, $H^p(\bbC^+)$ is a Banach space, while if
$0<p<1$, $H^p(\bbC^+)$ is only a quasi-Banach space.

Elements of $H^p(\bbC^+)$ possess 
non-tangential boundary limits almost everywhere on the real line
$\bbR$, and the space of their boundary values can
be identified with a closed subspace of $L^p(\bbR)$. 
More precisely, for $1<p<\infty$, the following decomposition holds
$$
L^p(\bbR) = H^p(\bbC^+)|_{\bbR}  \oplus
\ov{H^p(\bbC^+)|_{\bbR}}, 
$$
where $H^p(\bbC^+)\vert_{\bbR}$ denotes the boundary
values of $H^p(\bbC^+)$ and
$\ov{H^p(\bbC^+)\vert_\bbR}$ denotes its complex
conjugate.

The case $p=1$ is more delicate.  The {\em real variable} Hardy
space on $\bbR$ is defined as
$$
H^1(\bbR) = \big\{ f\in L^1(\bbR)\colon \, Hf\in  L^1(\bbR)\big\},
$$
where $H$ denotes the Hilbert transform, with norm $\|f\|_{H^1(\bbR)}=\|f\|_{L^1}+\|Hf\|_{L^1}$.
It turns out that 
$$
H^1(\bbR)=H^1(\bbC^+)\vert_{\bbR}\oplus\ov{H^1(\bbC^+)\vert_\bbR},
$$
where $H^1(\bbC^+)\vert_{\bbR},
\ov{H^1(\bbC^+)\vert_\bbR}$ isomorphically embed into $H^1(\bbR)$.
The space $H^1(\bbR)$ admits an atomic decomposition.  Precisely, every $f\in H^1(\bbR)$ 
can be written as $\sum_j\lambda_j a_j$, where $\lambda_j \in
\bbC$ with $\sum_j|\lambda_j|<\infty$ and $a_j$'s are
$H^1$-atoms, that is, they are $L^1$-functions which satisfy the
following three conditions: 
\begin{itemize}
\tbi there exists an interval $I_j$ such that $\text{supp}(a_j)\subseteq I_j$;
\tbi $\int_{I_j}a(x)dx=0$;
\tbi $|a_j(x)|\leq 1/|I_j|$, where $|E|$ denotes the Lebesgue measure of a measurable set
$E$.
\end{itemize}

Notably, every function $f \in H^1(\bbR)$ satisfies 
\begin{equation}\label{zero mean}
    \int_{\bbR}f(x)dx=0.
\end{equation}
\medskip

Let 
$\cS$  and
$\cS'$ denote the space of Schwartz functions and
the space of tempered
distributions on the real line $\bbR$ respectively.
For $f\in \cS'$ we equivalently denote by $\widehat f$ or $\cF f$ its Fourier transform, which, for $f\in\cS$, is given by
\[   \cF f(\xi)=\widehat f(\xi)
   =  \frac{1}{\sqrt{2\pi}} \int_{\bbR} f (x)e^{-ix\xi}\, dx. 
\]
The Fourier transform $\cF$ is an isomorphism of $\cS$
onto itself with inverse given by
$$
\cF^{-1}g(x) =
\frac{1}{\sqrt{2\pi}}\int_{\bbR} g(\xi)e^{ix\xi}\, d\xi.
$$
\medskip

The Hardy spaces can also be characterized using the Cauchy-Szeg\"o
projection,  If $f \in \cS$, the orthogonal projections
$P_\pm$ of $L^2(\bbR)$ onto $H^2(\bbC^+)$, $\ov{H^2(\bbC^+)}$,
resp., are given by 
\[
P_+f(x)=\lim_{y \to 0^+} \frac{1}{2\pi i}\int_{\bbR} \frac{f(t)}{t-x-iy}dt
\]
and
\[
P_-f(x)=\lim_{y \to 0^+} \frac{1}{2\pi i}\int_{\bbR} \frac{f(t)}{x-iy-t}dt.
\]
Note that 
$P_++P_-=\operatorname{Id}$.  By Riesz' theorem, for $p\in(1,\infty)$, the projections
$P_\pm$ extend to bounded linear operators of $L^p(\bbR)$ into
themself.  The 
Hilbert transform $H\colon L^2(\bbR)\to
L^2(\bbR)$ is defined by the identity
\begin{equation}\label{Hilb-tran:eq}
H=-i(P_+-P_-).
\end{equation}
By the boundedness of $P_\pm$, the Hilbert transform $H\colon L^p(\bbR)\to
L^p(\bbR)$ is also bounded when $p\in(1,\infty)$.

It is important to recall that, while $P_\pm$ fail to be bounded on $L^1(\bbR)$, 
$P_\pm \colon H^1(\bbR)\to H^1(\bbC^\pm)$ are
bounded, and that
\[
P_\pm (H^1(\bbR))=H^1(\bbC^\pm),     
\]
see \cite[Chapter 9]{zhu2007operator}.  Here and in what follows, we
use the natural identification of $\ov{H^2(\bbC^+)}$ with $H^2(\bbC^-)$.
\medskip

By
$\BMO(\bbR)$, or simply $\BMO$,  we denote the quotient space of functions  of bounded
mean oscillation modulo constants:
$$
\BMO(\bbR) = \Big\{ \vp \in L^1_{\rm loc}(\bbR):\, \|\vp\|_{\BMO(\bbR)}\colon= \sup_{I\subset\bbR}
\frac{1}{|I|} \int_I |\vp -\vp_I|\, dm <\infty \
\Big\},
$$
where $I\subset\bbR$ is any bounded interval,
$\vp_I=\frac{1}{|I|}\int_I \vp$ is the average of $\vp$ over
$I$. We emphasize that $\vp \in \BMO$ represents an equivalence class and
$\vp=\psi$ if and only if $\vp-\psi$ is constant.
Every function $f \in \BMO$ is Poisson integrable and, therefore, also belongs to $\cS'$.
We define the space of analytic $\BMO$
$$
\BMOA = \big\{ \vp\in \BMO(\bbR):\, \text{supp} (\widehat
\vp)\subseteq[0,\infty)\big\}.
$$
 It is well known that $\BMO(\bbR)$ is isomorphic to the dual of
$H^1(\bbR)$: for every bounded linear functional $\Lambda$ on $H^1(\bbR)$, there exists a unique function $\psi \in \BMO(\bbR)$ such that
$$
\Lambda(f)=\left\langle f,\psi\right\rangle\colon=\int_\bbR
  f(x)\ov{\psi(x)}\, dx, 
$$
for every $f \in \cS$, \cite[Chapter VI]{garnett2006}. Analogously, $\BMOA$ represents the dual of $H^1(\bbC^+)_{|_\bbR}$. In this article, we endow
$\BMO$ and $\BMOA$ with the dual norm.
\bigskip

This article is devoted to the study of the Toeplitz operators on
$H^1(\bbC^+)$
having bounded anti-analytic symbol. 
As pointed out by Janson, Peetre and Semmes in
\cite[Remark 4.1]{Janson1984}, the Toeplitz operator $T_{\ov\Theta}$ is
never bounded on $H^1(\bbC^+)$. 
\begin{prop}[Janson, Peetre and Semmes, Remark 4.1 \cite{Janson1984}]\label{Janson}
Let $\Theta \in H^\infty(\bbC^+)$ be non-constant. The Toeplitz
operator $T_{\ov\Theta}$  is unbounded on $H^1(\bbC^+)$.   
\end{prop}

While it is not difficult to check that the calculation indicated
in \cite{Janson1984} proves Proposition \ref{Janson}, we are going to
provide a different argument.
Indeed,  in \eqref{Formal adjoint}, we observe that the
formal adjoint of $T_{\ov\Theta}$ on $H^1(\bbC^+)$ is
the multiplication operator
$M_\Theta$ acting on $\BMOA$. Thus,  it becomes apparent that
$M_{\Theta}$ is not even well-defined as an operator from  $\BMOA$
into itself, as, in particular it fails to
map the equivalence class of $0$ into itself.

Therefore, we propose that the appropriate space for studying the
Toeplitz operator $T_{\ov\Theta}$ is the pre-dual of
$\BMOA/\operatorname{span} \{\Theta\}$. In this setting, $M_{\Theta}$
is well-defined as a map from $\BMOA$ to $\BMOA/\operatorname{span}
\{\Theta\}$. 
 
\begin{definition}\label{defn spazio H1theta} {\rm
 Let $\Theta\in H^\infty(\bbC^+)$. The space $H^1_\Theta$ is defined as the kernel of the linear functional induced by
$\Theta$ on $H^1(\bbC^+)$, that is, 
$$
H^1_\Theta\colon=\left\lbrace f\in  H^1(\bbC^+)\colon \left\langle f , \Theta\right\rangle = 0\right\rbrace.
$$
Note that $H^1_\Theta$ is a closed subspace of $H^1(\bbC^+)$, so
we endow $H^1_\Theta$ with its norm,
$$
\|f\|_{H^1_\Theta}=\|f\|_{H^1}.
$$
}
\end{definition}

It is clear that the dual of $H^1_\Theta$ is isomorphic to the
quotient space $\BMOA/\operatorname{span} \{\Theta\}$, see
\cite[Chapter 4]{rudin1991}. We note that when $\Theta$ is
constant, then $H^1_\Theta=H^1(\bbC^+)$, because of \eqref{zero
  mean}. \ms

In order to avoid trivialities, from now on we {\em always} assume
that the function $\Theta\in H^\infty$ we consider is {\em
  non-constant}. \ms

An
$H^\infty(\bbC^+)$-function $\Theta$ is called an {\em inner function} if
$|\Theta|=1$ a.e.\ on $\bbR$.  
Given an inner function $\Theta$, and $p\in[1,\infty)$, we  consider
the {\em Beurling subspace} 
\begin{equation}\label{Beur-sp:eq}
\Theta H^p\colon=\left\lbrace f \in H^p(\bbC^+) \colon
  (\ov\Theta f)_{|_\bbR} \in
  H^p(\bbC^+)\vert_{\bbR}\right\rbrace ,
\end{equation}
and the {\em backward-shift invariant subspace}
\begin{equation}\label{Beur-sp:eq2}
K^p_{\Theta}\colon=\left\lbrace f \in H^p(\bbC^+)
  \colon(\ov\Theta f)_{|_\bbR}  \in
  \ov{H^p(\bbC^+)}\vert_{\bbR}\right\rbrace . 
\end{equation}

The first main result provides some equivalent conditions to the
boundedness of $T_{\ov\Theta}$ on $H^1_\Theta$.

\begin{theorem}\label{equiv-cond:thm}
  Let $\Theta\in H^\infty(\bbC^+)$. 
  Then, the following conditions are equivalent.
  \begin{enumerate}[{\rm \quad (a)}]
\item  The Toeplitz operator
  $T_{\ov\Theta} \colon H^1_\Theta \to H^1(\bbC^+)$ is
  bounded;
  \item the Hankel operator $H_{\ov\Theta} \colon H^1_\Theta \to H^1(\bbC^-)$ is
  bounded;
  \item  the multiplication  operator $M_\Theta\colon \BMOA\to
    \BMOA/\operatorname{span} \{\Theta\}$ is bounded.
\end{enumerate}
If in addition $\Theta$ is an inner function, then {\rm (a) --
  (c)} are also equivalent to:
\begin{enumerate}[{\rm \quad (a)}]\setcounter{enumi}{3}
  \item we have the equality
      $H^1_{\Theta}=K^1_{\Theta}\oplus \Theta H^1$.
      \end{enumerate}
\end{theorem}

We point out that, if $p\in(1,\infty)$, it holds that $H^p(\bbC^+)
= K^p_{\Theta}\oplus  \Theta H^p$, see e.g.\ \cite{Nowak-Soltysiak}
for the case of the unit disk; the simple argument works also in the
case of $\bbC^+$. In Lemma \ref{H1Theta-prop:lem} we show that, when
$p=1$, $K^1_{\Theta}\oplus  \Theta H^1$ is a closed subspace of
$H^1_\Theta$, which of course, since we are assuming throughout the paper that $\Theta$ is
non-constant,
is properly contained in $H^1(\bbC^+)$.

 Particularly important inner functions are the meromorphic
inner functions, see \cite[Chapter VII]{levin1964}.
The second main result of this article shows that the
Toeplitz operator associated to the meromorphic inner function
$\Theta=e^{i\tau (\cdot)}$, $\tau>0$, is bounded.

\begin{theorem}\label{main-result}
The Toeplitz operator $T_{e^{-i\tau (\cdot)}}\colon H^1_{e^{i\tau(\cdot) }}
\to
H^1(\bbC^+)$ is bounded.
\end{theorem}

The rest of the paper is organized as follows.  In Section \ref{2:sec}
we discuss some general properties of the Toeplitz operators on
$H^1(\bbC^+)$ with anti-analytic symbol, 
 we describe some properties of the subspace
$H^1_\Theta$, and we conclude with  proving Theorem  \ref{equiv-cond:thm}.  Section
\ref{3:sec} is devoted to the proof of Theorem \ref{main-result}.  We
conclude in 
Section \ref{open-question:sec} with some final remarks and 
 open questions.
\ms

Our interest in the Toeplitz operators stems from our work in
\cite{Bellavita2024}, where we characterized the dual of the
$1$-Bernstein spaces (also called Paley--Wiener 1-spaces).
Precisely, we have used Theorem
\ref{main-result} in \cite[Theorem 1]{Bellavita2024} to provide an
explicit description of the 
dual of $1$-Bernstein spaces in terms of $\BMO$ entire functions.
 We refer the reader to
\cite[Lectures 20-21]{Levin-Lectures},
\cite{Plancherel-Polya2,Plancherel-Polya1}, e.g., 
 for the properties of the Bernstein spaces. Further results can be
 found in \cite{MPS}.

The $1$-Bernstein spaces are particular examples of $1$-de Branges
 spaces.  See \cite{Baranov2006} and \cite{Bellavita20} for some properties of the $1$-de Branges
 spaces.  In 
the upcoming paper \cite{Bellavita2025}, we apply this circle of ideas to provide a
characterization of the dual of the 1-de Branges spaces.

In \cite{Aleksandrov} A.\ B.\ Aleksandrov characterized the
backward-shift invariant subspaces of $H^1(\bbD)$, we refer to the
monograph \cite{Cima-Ross} for a detailed proof.
We also point out that some descriptions and properties of the dual
 of the  $1$-backward-shift invariant  subspaces in the unit disk, which are
 strictly related to the $1$-de Branges spaces, have been provided in
 \cite{OLOUGHLIN2022}, \cite{BESSONOV201562}  and
 \cite{Dyakonov2022}. 
\medskip  

\section{Toeplitz operators on $H^1(\bbC^+)$ with anti-analytic symbol}\label{2:sec}

In this section we present the Toeplitz operators of
$H^1(\bbC^+)$ and we focus our attention on those with
anti-analytic symbol.
However, before proceeding with the case of the upper half-plane, we
briefly recall the situation in the case of the unit disk.
\ms

Let $\cP$ denote the space of trigonometric polynomials on $\bbT$, and
$\cP'$ its dual. Then, $\cP'$ is the space of formal Fourier series. 
For every $f\in \cP$,  we consider the
multiplication operator
$$
M_b\colon \cP\ni f \mapsto bf \in \cP' .
$$ 
Since the action of the Szeg\"o projection $P_+$ is well defined on
$\cP'$, we define $T_b\colon= P_+M_b$, and observe that
\begin{equation} 
\label{Toeplitz operator disk}
T_b\colon \cP\to \cP'
\end{equation}
is well defined.
We now introduce the $\log$-BMO space, that we denote by $\LMO$.
An $L^1(\bbT)$ function $f$ belongs to $\LMO$ if
\begin{equation*}\label{E:stegenga 2}
\sup_{I\subseteq\bbT} \frac{\log(1/|I|)}{|I|}\int_{I}|f(e^{i\theta})-f_I|\frac{d\theta}{2\pi}<\infty,
\end{equation*}
where $I$ ranges over all subarcs of $\bbT$. D.\ Stegenga initially
characterized bounded Toeplitz operators on $H^1(\bbT)$ under certain
additional assumptions on the symbol $b$. Subsequently, Janson,
Peetre, and Semmes provided a complete characterization. 

\begin{theorem}[Stegenga, Theorem 1.1 \cite{Stegenga1976} \& Janson, Peetre and Semmes, Theorem 2.i \cite{Janson1984}]
Let 
$$
b(e^{i\theta})=\sum_{n\in\bbZ}\hat{b}(n)e^{in\theta}\in \cP'. 
$$
Then, the Toeplitz operator $T_b$ defined as in \eqref{Toeplitz operator
  disk} admits a bounded extension as an operator
$T_b\colon H^1(\bbT)\to H^1(\bbT)$ if and only if
$b \in L^\infty(\bbT)\cap \LMO$.
\end{theorem}

Next, if $\vp \in H^\infty$, it is easy to see that the adjoint of
$T_{\ov{\vp}} \colon H^2(\bbT)\to H^2(\bbT)$
is the multiplication operator $M_{\vp}$. 
Therefore, the boundedness of $T_{\ov{\vp}}$ on
$H^1(\bbT)$ implies the boundedness of $M_\vp$ on
$\BMOA(\bbD)$.
Actually, the reverse implication is true as well.
\begin{theorem}[Stegenga, Theorems 1.2 and 3.8 \cite{Stegenga1976}] The  linear operator
$M_\vp$ is a bounded multiplier of
$\operatorname{BMOA}(\bbD)$ if and only if $\vp \in H^\infty$ and its
boundary values belong to $\LMO$.     In particular, if $\Theta$ is an
inner function, $T_{\ov\Theta}$ is bounded if and only if $\Theta$ is
a finite Blaschke product. 
\end{theorem}
\ms

We now turn to the case of  the upper half-plane. Our first step is to
identify  the class of
allowable symbols of the Toeplitz operators.  Notice that, in the case
of the unit circle, $P_+\colon\cP\to \cP$, so that we could
consider symbols in $\cP'$.  In the case of the real line,
$P_+\colon\cS\not\to \cS$ and the situation is a bit more
involved. 

Let $\vp\in L^1_{\text{loc}}(\bbR)$ be a moderate increasing
function, that is, at most of polynomial growth.  
For every $f \in \cS$, we consider the
multiplication operator
$$
M_\vp\colon \cS\ni f \mapsto \vp f\in L^2(\bbR).
$$
Next, we observe that the space
$$
\cS_+= \big\{ f=P_+\phi\colon \phi \in\cS,\, \supp(\widehat\phi\,)\subseteq
[0,\infty)\big\}
$$
is contained and dense in $H^p(\bbC^+)$, for all $1\le p<\infty$.
Notice in particular that $\cS_+$ is a space of holomorphic functions
in $\bbC^+$, given by the Cauchy--Szeg\"o integral of Schwartz functions
whose Fourier transform is supported on the half-line $[0,\infty)$. Therefore, 
if $f=P_+\phi$, then $f_{|_\bbR}=\phi$. 

For  $\vp$ in $L^1_{\text{loc}}(\bbR)$ at most of
moderate growth, we define
\begin{equation}\label{e:toeplitz semipiano}
 T_\vp\colon=P_+M_\vp \colon \cS_+ \to \Hol(\bbC^+).
\end{equation}
Analogously, we consider the {\em Hankel operator with symbol}
$\vp$ 
\begin{equation}\label{Hankel-op-C+:eq}
 H_\vp\colon=P_- M_\vp \colon \cS_+ \to \Hol(\bbC^-).
\end{equation}

 We are interested in studying under which conditions on the symbol of the
Toeplitz operator,  $T_\vp$ extends to a bounded linear operator
$T_\vp \colon H^1(\bbC^+)\to H^1(\bbC^+)$.

Clearly, if $\vp=u \in H^\infty(\bbC^+)$, then
 $T_u$ is bounded on $H^1(\bbC^+)$. 
On the other hand, by
Proposition~\ref{Janson} we know that if
$\vp$ is an anti-analytic symbol,  
$\vp=\ov\Theta$, with
$\Theta \in H^\infty(\bbC^+)$, then $T_{\ov\Theta}$ is
unbounded. However, this fact can be also observed directly.  Indeed,
if 
$T_{\ov\Theta}\colon H^1(\bbC^+)\to H^1(\bbC^+)$ were bounded,
then 
$T_{\ov\Theta}^* \colon \BMOA\to\BMOA$ would also be bounded.  However, it is
simple to see that $T_{\ov\Theta}^*=M_\Theta$.  Indeed, 
if $f,g \in \cS_+$, we have that
\begin{equation}\label{Formal adjoint}
\la T_{\ov\Theta}f,g\ra = \la \ov\Theta f,g \ra =\la  f,M_\Theta g \ra.
\end{equation}
However, the multiplication by
$\Theta$ is not a well-defined operator on $\BMOA(\bbC^+)$. Recall that
elements of $\BMOA(\bbC^+)$ are equivalent classes of
functions modulo constants. Then,  $M_{\Theta}$ maps elements 
of $\BMOA(\bbC^+)$ into functions modulo multiples of
$\Theta$. Therefore, 
the natural question is for which $\Theta\in H^\infty$, or more
generally even $\Theta\in\BMOA$, 
$$
M_\Theta\colon \BMOA\to \BMOA/\operatorname{span}\{\Theta\}
$$
is bounded. Note that $\BMOA/\operatorname{span}\{\Theta\}$ is the dual
of the annihilator of $\Theta$ in $H^1(\bbC^+)$, that is the space
$H^1_\Theta$, see Definition~\ref{defn spazio H1theta}.
Note also that this situation does not occur in the case of $\bbD$.
Indeed, $\BMO$ is {\em not} a quotient space, and the multiplication
operator on $\BMO$ produces a well defined function.

The above
discussion is the motivation for the next definition.  First a lemma.
\begin{Lemm}\label{H1ovTheta-domain:lem}
  Let $\Theta\in H^\infty$. Then,
  $H^1_\Theta\cap \cS_+$ is dense in $H^1_\Theta$.
\end{Lemm}

\proof
 Let $f\in H^1_\Theta$, and let $\{h_n\}\subseteq\cS_+$ be such
that $h_n\to f$ in $H^1(\bbC^+)$.  Then, $\la h_n,\Theta\ra \to 0$ as
$n\to\infty$.  Let $h_0\in\cS_+$ be such that $\int_\bbR h_0\ov\Theta=1$.  Then,
$g_n:=h_n -\la h_n,\Theta\ra h_0\in\cS_+\cap H^1_\Theta$, and $g_n\to
f$ in $H^1(\bbC^+)$, as $\to\infty$. 
This proves the lemma.
\epf

\begin{definition}{\rm Let $\Theta\in H^\infty$.  Then we
    define the domain of the
    Toeplitz operator $T_{\ov\Theta}$ as $H^1_\Theta\cap \cS_+$.
  }
  \end{definition}
  In other words, given an  inner function $\Theta$, we always think
  of the Toeplitz operator  $T_{\ov\Theta}$ as initially defined on
  $H^1_\Theta\cap \cS_+$
  which is a dense
  subspace of 
  $H^1_\Theta$, where the latter is endowed with its topology as closed
  subspace  of $H^1(\bbC^+)$. We
  wish to determine for 
  which $\Theta\in H^\infty$,
$T_{\ov\Theta}$ extends to a bounded linear operator 
$T_{\ov\Theta}\colon H^1_\Theta\to H^1(\bbC^+)$.
With Theorem~\ref{main-result} that we prove below, we show
that the answer is in the positive when $\Theta=e^{i\tau(\cdot)}$.
Analogously, we think of the Hankel operator $H_{\ov\Theta}$ as defined
on the dense subspace of $H^1_\Theta$, $H^1_\Theta\cap \cS_+$.  
\medskip

We provide some conditions equivalent to the boundedness of
$T_{\ov\Theta}$, when $\Theta$ is an inner function.
Recall that the subspaces $\Theta H^1$ and $K^1_{\Theta}$ are defined
in \eqref{Beur-sp:eq} and \eqref{Beur-sp:eq2}, resp.

\begin{Lemm}\label{H1Theta-prop:lem}
  Let $\Theta$ be an inner function. Then, we have
  \begin{enumerate}[{\rm \quad (a)}]
\item $ K^1_\Theta\oplus\Theta H^1$ is a closed subspace of
  $H^1_\Theta$;
    \item $K^1_\Theta\subseteq \ker(T_{\ov\Theta})$.
    \end{enumerate}
\end{Lemm}

\begin{proof}
It is well known, and immediate to check, that $\Theta H^1$ and
$K^1_\Theta$ are closed subspaces of $H^1(\bbC^+)$.
If $f\in K^1_\Theta\cap\Theta H^1$,
then there exists $g,h\in H^1(\bbC^+)$ such that
$f_{|_\bbR}= (\Theta \ov g)_{|_\bbR}$ and $f_{|_\bbR}= (\Theta
h)_{|_\bbR}$. This implies that $\ov g_{|_\bbR}=h_{|_\bbR}$, whence
$g=h=0$.  Therefore, $ K^1_\Theta\oplus\Theta H^1$ is a direct sum of
closed subspaces of $H^1(\bbC^+)$. However, if $f \in K^1_\Theta$,
$f_{|_\bbR}= (\Theta \ov g)_{|_\bbR}$ for a $g\in H^1(\bbC^+)$, then
$$
\la f,\Theta\ra = \int_\bbR f\ov \Theta = \int_\bbR \ov{g} =0,
$$
so that $f\in H^1_\Theta$.  Analogously, if $f \in \Theta H^1$, and $f_{|_\bbR}= (\Theta
h)_{|_\bbR}$, with $h\in H^1(\bbC^+)$, then
$$
\la f,\Theta\ra = \int_\bbR h = 0.
$$
Thus, $ K^1_\Theta\oplus\Theta H^1$ is a closed subspace of
$H^1_\Theta$, and (a) follows.

(b)
Notice that Lemma \ref{H1ovTheta-domain:lem} guarantees that $T_{\ov\Theta}$ is defined on a
dense subspace of $H^1_\Theta$. 
Let $f\in K^1_\Theta$, so that $f_{|_\bbR}= (\Theta \ov
g)_{|_\bbR}$ for some $g\in H^1(\bbC^+)$. Then,
$$
T_{\ov\Theta}f = P_+(\ov g) = 0.
$$
This proves the lemma.
\end{proof}

\begin{proof}[Proof of Theorem~\ref{equiv-cond:thm}] 
In order to prove that (a) and (b) are equivalent, note that
\begin{align*}
  \|{T_{\ov\Theta}f}\|_{H^1(\bbC^+)}
  & = \|(\text{Id}-P_-)\left( \ov\Theta f \right)\|_{H^1(\bbC^+)}
    \leq
    \|\Theta\|_{H^\infty}\|f\|_{{H^1(\bbC^+)}}+\|H_{\ov\Theta}(f)\|_{{H^1(\bbC^-)}}    .
\end{align*}
Hence (b) implies (a). 
For the reverse implication, note that $P_- = \text{Id}-P_+$, so that
\begin{align*}
  \|{H_{\ov\Theta}f}\|_{H^1(\bbC^-)}
  & = \|(\text{Id}-P_+)\left( \ov\Theta f \right)\|_{H^1(\bbC^-)}
    \leq
    \|\Theta\|_{H^\infty}\|f\|_{{H^1(\bbC^+)}}+\|T_{\ov\Theta}(f)\|_{{H^1(\bbC^+)}}    .
\end{align*}

Next we prove that (a) and (c) are equivalent.
  If $T_{\ov\Theta}$ is bounded, then $M_\Theta$ is bounded
  and $\|M_\Theta\|=
\|T_{\ov\Theta}\|$.
Conversely, assume that $M_{\Theta}$ is bounded. Consider, $g \in  H^1_\Theta$. Then,
\begin{align*}
  \|T_{\ov\Theta}(g)\|_{H^1(\bbC^+)}
  & =\sup_{ v \in\VMOA,\, \|v\|=1}|\left\langle v ,  T_{\ov\Theta}(g) \right\rangle|\\
&=\sup_{ v \in\VMOA,\, \|v\|=1}|\left\langle \Theta v ,   g\right\rangle| \leq \|M_{\Theta}\| \|g\|_{H^1_\Theta},
\end{align*}
where $\VMOA$ is the pre-dual  of $H^1(\bbC^+)$.
Therefore,
$T_{\ov\Theta}$ is bounded and, by the first part of the proof, $\|T_{\ov\Theta}\|=
\|M_{\Theta}\| $. Hence, (a) and (c) are equivalent.
\ms

Suppose now that $\Theta$ is an inner function and that (d) holds. 
Let $f\in H^1_\Theta$, 
$f=f_1+f_2 \in K^1_\Theta\oplus \Theta H^1$, where $f_1 \in
K^1_\Theta$, $f_2 \in  \Theta H^1$ and $f_2= \Theta g$, with $g\in
H^1(\bbC^+)$. 
By Lemma \ref{H1Theta-prop:lem} (b), 
we have that  
$$
T_{\ov\Theta} f= T_{\ov\Theta} f_2= g \in H^1(\bbC^+)
$$
and
$$
\|T_{\ov\Theta} f\|_{H^1} = \| g\|_{L^1(\bbR)} 
=\|\ov \Theta f_2\|_{L^1(\bbR)} =\|f_2\|_{L^1(\bbR)}  =
\|f_2\|_{H^1(\bbC^+)} \leq \|f\|_{H^1_\Theta}.
$$
Hence, (d) implies (a).   Finally, suppose (a) holds.  Consider $f
\in H^1(\bbC^+)$.  We have that
$\ov\Theta
f\in H^1(\bbR)$. 
Indeed, clearly $\ov\Theta
f\in L^1(\bbR)$, while 
$$
H(\ov\Theta f) = -i\big(2P_+ -\operatorname{Id}\big)(\ov\Theta f) =
-2i T_\Theta f +i \ov\Theta f \in L^1(\bbR),
$$
since (a) holds.
Since $\ov\Theta
f\in H^1(\bbR)$, it follows that
\begin{equation}\label{above:eq}
\ov\Theta
f =  P_- (\ov\Theta f) +P_+(\ov\Theta f),
\end{equation}
that is,
$$
f=  \Theta P_- (\ov\Theta f) +\Theta P_+(\ov\Theta f) =\colon f_1+f_2,
$$
where $f_2 =\Theta T_{\ov\Theta} f \in \Theta H^1\subseteq H^1$.  
Therefore, 
$$
\| f_2\|_{H^1} \le \|T_{\ov\Theta} f \|_{H^1} \le C \| f \|_{H^1}.
$$
Moreover, since
 $f_1=f-f_2$, $f_1 \in H^1(\bbC^+)$,
$$ 
f_1 =\Theta P_- (\ov\Theta f) \in K^1_\Theta
$$
and,  we also have 
$$
\| f_1\|_{H^1} \le (1+C) \| f \|_{H^1}.
$$
This shows  that the maps $H^1_\Theta \ni f \mapsto f_1\in K^1_\Theta$ and 
$H^1_\Theta \ni f \mapsto f_2\in \Theta H^1$ are bounded,  and that $H^1_\Theta = K^1_\Theta+\Theta
H^1$.   The
conclusion now follows from Lemma \ref{H1Theta-prop:lem} (a).  
Notice however that,  since
$\ker(T_{\ov\Theta})\cap \Theta H^1 =\{0\}$, while
$K^1_\Theta\subseteq \ker(T_{\ov\Theta})$,  we also have that
$$
K^1_\Theta=\ker(T_{\ov\Theta}).
$$
This completes the proof.
\end{proof}
\ms

\section{Boundedness of $T_{e^{-i\tau (\cdot)}}$ on $H^1_{e^{i\tau(\cdot) }}$}\label{3:sec}

Before proving Theorem \ref{main-result}, we need some preliminary considerations. 
If $J$ is  a subset of $[0,\infty)$, we denote by $H^1_J$  the
subspace of $H^1(\bbC^+)$ of functions whose Fourier transform
of its boundary values vanishes on $J$.
Observe that if $J$ is a singleton, $J=\{\tau\}$, with $\tau>0$,
then $H^1_{\{\tau\}}= H^1_{e^{i\tau(\cdot) }}$. Therefore, from now
on, we use the former notation.
\begin{Lemm}\label{Claim1}
Let
$$
\cS_{0,\tau} \colon=\big\{ \vp\in\cS:\
\supp(\widehat\vp)\subseteq[0,\infty),\, \textstyle{\int \vp(x)x^k dx= \int
  \vp(x)  x^j e^{-i\tau x} dx =0} \ \forall j,k\in\bbN\big\}.
$$
Then, $\cS_{0,\tau} $ is dense in $H_{\{\tau\}}^1$. 
\end{Lemm}
\begin{proof}
In order to prove the lemma, we begin by showing that the subspace
$$
H^1_{00} =\big\{\phi\in L^1\colon \widehat\phi\in C^\infty_c(\bbR),\,
\operatorname{dist}(\supp(\widehat\phi),\{0\})>0
  \big\} 
  $$
  is dense in $H^1(\bbR)$.  This is \cite[Corollary
  pag.~230]{stein1970}. 
The proof is based on the fact that a function in $L^1_0 =\{g\in
L^1(\bbR)\colon \widehat g(0)=0\}$ can be approximated by $L^1$-functions
whose Fourier transform is compactly supported, with support having
positive distance from the origin.   We repeat the proof here for the
reader's convenience, and since we use a similar argument in the
second step of the proof.

As customary, for 
$r>0$ and a function $g$, we denote by $g_r$ and $g^r$ the functions
defined by  
$$
g_r(\xi) = \frac1r g(\xi/r) \text{ and } g^r(\xi) =g(r\xi) .
$$
Recall that $\widehat{g_r \,}={\widehat g\,}^r$.
Let $\eta\in C^\infty_c$, $\eta(\xi) =1$
if $|\xi|\le 1$ and $\eta(\xi) =0$ if $|\xi|\ge2$. For $r>0$ we define
the operator  $T_r f= \cF^{-1}\big(\eta^{1/r}\big)*f$. 
Notice that  
$$
\|T_r f\|_{L^1} \leq \|\big(\cF^{-1}\eta\big)_{1/r}\|_{L^1} \| f\|_{L^1} = \|\cF^{-1}\eta\|_{L^1} \|f\|_{L^1}.
$$
Hence,
the operators $T_r$ have uniformly bounded norm on $L^1$.
An elementary argument
shows that,
if $f\in L^1_0$, 
$T_\epsilon f \to 0$ as $\epsilon \to 0$.  Therefore, if $f\in L^1_0$,
$$
T_R(\text{Id}-T_\epsilon)f\to f \qquad\text{as } R\to\infty,\, \epsilon\to 0.
$$

Moreover, 
$$
\cF\big(T_R(\text{Id}-T_\epsilon)f\big) =
\eta^{1/R}(1-\eta^{1/\epsilon})\widehat f
$$
has support  contained in
$\{\xi:\, \epsilon \leq |\xi|\le 2R\}$.
Let $\widetilde \eta$ be another smooth cut-off function, having
support contained in $\{|\xi|\le \epsilon'\}$.
Then, the functions 
$$
f_{R,\epsilon}\colon=T_R(\text{Id}-T_\epsilon)f \in L^1_0
$$ 
and for
$0<\epsilon'<\epsilon$, 
 $\widetilde \eta_{\epsilon'}*\widehat{f_{R,\epsilon}} \in C^\infty_c$  with support
contained in 
$$
\{\xi:\, (\epsilon-\epsilon')\leq |\xi|\leq 2R+\epsilon'\}.
$$ 
Then,
$H^1_{00}$ is dense in $L^1_0$. However, $H^1_{00}$ is dense also in
$H^1(\bbR)$. For, given any function $g\in
H^1_{00}$, let $\chi\in C^\infty_c(\bbR)$, $\chi=-i\, \text{sign}$ on
$\text{supp}(\widehat g)$. Then,  
$$
Hg = \cF^{-1}\big({\chi}\widehat{g}\big)\in L^1(\bbR)
$$ (recall, $H$ denotes the Hilbert transform).

Now, given  $f\in H^1(\bbR)$, let $f_n\colon= \widetilde
\eta_{2/n}*f_{n,1/n}$, with the above notation. Then, $f_n\to f$, and $Hf_n = \widetilde \eta_{2/n}* T_n(\text{Id}-T_{1/n})Hf\to Hf$ in
$L^1$, as
$n\to\infty$. 
(This is the argument  taken from \cite{stein1970}.)
\medskip

For $\tau\in\bbR$, consider
$$
H^1_{\tau 0} =\big\{\phi\in L^1\colon \widehat\phi\in C^\infty_c(\bbR),\,
\operatorname{dist}(\text{supp}(\widehat\phi),\{\tau\})>0
  \big\} .
  $$
  Notice that
$H^1_{\tau 0} = e^{i\tau(\cdot)}H^1_{00}$. Also we set
$H^1_{\{0,\tau\},0} \colon= H^1_{00}*H^1_{\tau 0}$. Then, 
  $$
H^1_{\{0,\tau\},0 }=\big\{\phi\in L^1\colon \widehat\phi\in C^\infty_c(\bbR),\,
\operatorname{dist}(\text{supp}(\widehat\phi),\{0, \tau\})>0
  \big\} .
  $$
Indeed, it is clear that $H^1_{0
    0}*H^1_{\tau 0} \subseteq\{\phi\in L^1\colon \widehat\phi\in C^\infty_c(\bbR),\,
\operatorname{dist}(\text{supp}(\widehat\phi),\{0, \tau\})>0\}$.
On the other hand, if $\eta\in  C^\infty_c(\bbR)$, and
$\operatorname{dist}(\text{supp}(\eta),\{0, \tau\})>0$, then we can find
$v_1,v_2 \in  C^\infty_c(\bbR)$, $v_1,v_2=1$ on $\text{supp}(\eta)$ and
$v_1$ identically $0$ in a neighborhood of the origin, and $v_2$
identically $0$ in a neighborhood of $\tau$. Then, $\eta=v_1(v_2\eta)$,
Moreover, $\cF^{-1}v_1 \in H^1_{00}$, $\cF^{-1}(v_2\eta)
\in H^1_{\tau 0}$. and $\cF^{-1}\eta = (\cF^{-1}v_1)*(\cF^{-1}(v_2\eta))$.

Finally, in order to prove the lemma, it suffices to show that 
$H^1_{\{0,\tau\},0 }$ is dense in $H^1_{\{\tau\}}$. Let $f\in
H^1_{\{\tau\}}$. In particular, $f, e^{-i\tau(\cdot)}f\in L^1_0$ so that
$$
T_n(\text{Id}-T_{1/n})f\to f,\quad  e^{i\tau(\cdot)}T_n(\text{Id}-T_{1/n})(
e^{-i\tau(\cdot)}f)\to f
$$
as $n\to \infty$, which imply that 
\begin{align*}
T_n(\text{Id}-T_{1/n})[ e^{i\tau(\cdot)}T_n(\text{Id}-T_{1/n})(
e^{-i\tau(\cdot)}f)] \to f.
\end{align*}
Once suitable regularized with $\widetilde\eta$ as before, the
functions on the left hand side above are in $H^1_{\{0,\tau\},0}$ and
the conclusion follows by the same argument as before.
This proves the lemma.
\end{proof}

\begin{prop}\label{Claim2}
It holds that
\begin{equation}\label{Claim1-v5} 
H_{\{\tau\}}^1 = H^1_{[\tau,\infty)} \oplus H^1_{[0,\tau]},
\end{equation}
where, if $H_{\{\tau\}}^1\ni f=f_1+f_2\in H^1_{[\tau,\infty)} \oplus H^1_{[0,\tau]} $, then
$$
\|f\|_{H^1}=\|f\|_{H_{\{\tau\}}^1}\asymp \|f_1\|_{H^1}+\|f_2\|_{H^1}.
$$
\end{prop}
\begin{proof}
Clearly,
$H^1_{[\tau,\infty)}$ and $H^1_{[0,\tau]}$ are closed subspaces of  $H_{\{\tau\}}^1$
and they have trivial intersection.
Thus, the sum in~\eqref{Claim1-v5} is a direct sum and it is a closed
subspace of $H_{\{\tau\}}^1$.  If it is  a proper subspace, there exists
$f\in H_{\{\tau\}}^1\setminus\big( H^1_{[\tau,\infty)} \oplus H^1_{[0,\tau]} \big)$.
Then, there exists $b\in \BMOA$ such that 
$b\vert_{H^1_{[\tau,\infty)} \oplus H^1_{[0,\tau]}}=0$ and $\left\langle f ,
b\right\rangle \neq 0$. 

Let 
$\{g_n\}\subseteq \cS_{0,2\pi}$, $g_n\to f_0$ in
$H^1(\bbR)$. 
Then, according to Lemma \ref{Claim1},  $\widehat g_n= \psi_{n,1}+\psi_{n,2}$, where
$\psi_{n,1},\psi_{n,2}\in\cS$, $\text{supp}(\psi_{n,1})\subseteq[0,\tau]$, and
$\text{supp}(\psi_{n,2})\subseteq[\tau,\infty)$. For $j=1,2$, we define
$$
F_{n,j}(z) =\frac{1}{2\pi} \int_0^\infty e^{iz\xi} \psi_{n,j}(\xi)\, d\xi
$$
and $F_n = F_{n,1}+F_{n,2}$, so that $F_{n,j}\in H^1(\bbC^+)$,
$ F_{n,1} \in H^1_{[\tau,\infty)}$, $ F_{n,2} \in H^1_{[0,\tau]}$, and 
$F_n\to f$ in $H^1(\bbC^+)$.  Therefore,
\begin{align*}
\left\langle f ,  b\right\rangle & =\lim_{n\to\infty} \left\langle F_{n,1}+F_{n,2} , b\right\rangle=
 \lim_{n\to\infty} \left\langle  F_{n,1},  b\right\rangle +\left\langle F_{n,2} , b\right\rangle
  = 0,
  \end{align*}
a contradiction.  The equivalence of the two norms follows from the open mapping theorem.
\end{proof}

The proof of Theorem \ref{main-result} is now immediate.

\proof[Proof of Theorem \ref{main-result}] Since, if $\Theta=e^{i\tau (\cdot)}$, 
$$
K^1_\Theta  =H^1_{[\tau,\infty)} \quad \text{ and}\quad  \Theta H^1 =H^1_{[0,\tau]},
$$
the conclusion follows from Theorem \ref{equiv-cond:thm} and
Proposition \ref{Claim2}.
\epf

\ms

\section{Open questions}\label{open-question:sec}

\begin{openproblem}{\rm 
  Clearly, the first question we are interested in is whether there are inner
functions $\Theta$, other than the meromorphic inner
functions $e^{i\tau (\cdot)}$,
such that the associated Toeplitz operator $T_{\ov\Theta}\colon
H^1_{\Theta}\to H^1(\bbC^+)$ is bounded, or equivalently, for which inner functions $\Theta$ the equality
$$ %
 H^1_{\Theta}=K^1_{\Theta}\oplus \Theta H^1
$$
holds.
Moreover, from Theorem \ref{equiv-cond:thm} we know that, if $\Theta$
is an inner function,
$T_{\ov\Theta}\colon
H^1_{\Theta}\to H^1(\bbC^+)$ is unbounded if and only if
$$
K^1_{\Theta}\oplus \Theta H^1 \subsetneq H^1_{\Theta},
$$
that is, if and only if the quotient $H^1_\Theta/(K^1_{\Theta}\oplus \Theta H^1 )$ is
non-trivial. When non-trivial, it would be interesting to describe it.

}
\end{openproblem}

\begin{openproblem}{\rm 
The second question delves deeper into the subject and it is related
to another interesting problem. 
In analogy with the holomorphic case, we introduce a subspace of the
real variable Hardy space $H^1(\bbR)$.   For $b\in\BMO(\bbR)$ we
consider the pre-annihilator of $\{b\}$, that is,
\begin{equation} \label{H1b:eq}
H^1_b(\bbR)\colon=\left\lbrace \phi\in  H^1(\bbR)\text{
    such that } \la \phi, b\ra = 0\right\rbrace.
\end{equation}
Observe that, if $b=\Theta\in H^\infty(\bbC^+)$
$$
H^1_{\Theta}(\bbR) =\ov{H^1(\bbC^+)}_{|_\bbR}\oplus (H^1_\Theta)_{|_\bbR}.
$$
Following the ideas in Theorem \ref{equiv-cond:thm}, we
reformulate the boundedness of $T_{\ov\Theta}$ in terms of
$[\ov\Theta,H]$.  

\begin{prop}
  Let $\Theta$ be an inner function. Then, the following are equivalent:
\begin{enumerate}[{\rm \quad (a)}]
\item  the Toeplitz operator
  $T_{\ov\Theta} \colon H^1_\Theta \to H^1(\bbC^+)$ is bounded;
  \item  the commutator
$[\ov\Theta,H]\colon H^1_{\Theta}(\bbR)\to L^1(\bbR)$
is bounded.
\item the multiplication operator $M_{\ov\Theta}\colon H^1_\Theta(\bbR)\to H^1(\bbR)$ is
  bounded.
\end{enumerate}
\end{prop}

\begin{proof}
  We note that for every $H^1_{\Theta}(\bbR)\ni f=f_1+f_2$, with
  $f_1 \in \ov{H^1(\bbC^+)}_{|_\bbR}$ and $f_2 \in
  H^1_\Theta $, by using \eqref{Hilb-tran:eq} we have that
\[
[\ov\Theta,H](f_1)= \ov\Theta H(f_1)
-H(\ov\Theta f_1) = \ov\Theta(-if_1)+i\ov\Theta f_1=0
.
\]
Moreover, 
\[
[\ov\Theta,H](f_2)= \ov\Theta H(f_2) -H(\ov\Theta f_2)= -i\ov\Theta f_2+iP_+(\ov\Theta f_2)-iP_-(\ov\Theta f_2)=-2iH_{\ov\Theta}(f_2).
\]
Consequently, $[\ov\Theta,H]$ is bounded if and only if the
Hankel operator $H_{\ov\Theta}$ is bounded from $H^1_{\Theta}$
to $\ov{H^1(\bbC^+)\vert_{\bbR}}$, which is
equivalent to the boundedness of 
$T_{\ov\Theta}\colon H^1_\Theta\to H^1(\bbC^+)$, by
Theorem \ref{equiv-cond:thm}.  Thus, (a) and (b) are equivalent.

Next, suppose 
(a) holds.  Then we have that 
if
 $f\in H^1_\Theta$, then $(\ov\Theta f)_{|_\bbR}\in
 H^1(\bbR)$. Indeed, 
we wish to show that
$(\ov\Theta f)_{|_\bbR}$ and $H((\ov\Theta f)_{|_\bbR}) \in L^1(\bbR)$.  Clearly $\ov\Theta
f\in L^1(\bbR)$, and 
$$
H(\ov\Theta f) = -i\big(2P_+ -\operatorname{Id}\big)(\ov\Theta f) =
-2i T_\Theta f +i \ov\Theta f \in L^1(\bbR),
$$
since (a) holds.  
Therefore, $(\ov\Theta f)_{|_\bbR}\in H^1(\bbR)$, and
$$
\|\ov\Theta f \|_{H^1(\bbR)} \le 2 \|T_\Theta f \|_{H^1(\bbC^+)} +
\|f\|_{H^1(\bbC^+)} \le C \|f\|_{H^1(\bbC^+)},
$$
that is, (c) holds.  
The fact that (c) implies (a) follows from the fact that $P_+
\colon H^1(\bbR)\to H^1(\bbC^+)$, and the lemma follows.
\end{proof}

From the above Lemma, the following question arises:\\ 
For which inner function $\Theta$ is the commutator
$[\ov\Theta,H]\colon H^1_{\Theta}(\bbR)\to L^1(\bbR)$  bounded?    
}
\end{openproblem}

\begin{openproblem}{\rm
    In the real-variable theory of Hardy spaces, it is natural to
    extend Question 2 to more general symbols.  If $b\in\BMO(\bbR)$,
    then it has been known that the commutator $[b,H]$ is not bounded
    from $H^1(\bbR)$ to $L^1(\bbR)$.  It was then natural to ask on
    which subspace of   $H^1(\bbR)$ the 
commutator $[b,T]$ is bounded, where $T$ denotes a general Cald\'eron--Zygmund operator.
In \cite{perez1995}, C.\ Perez gave the following definition.
\begin{definition}{\rm 
    Let $b\in \BMO$ be given. A function $a$ is called a $b$-atom if
it is an atom of $H^1(\bbR)$ (see Section~\ref{Sec1})  with the additional
property that $\int_\bbR a\ov b = \la a,b\ra=0$.   
The space $H^1_{b,{\rm at}}$ is defined to be the space of
$L^1(\bbR)$-functions $f$ which can be written as
$$
f =\sum_{j}\lambda_j a_j ,
$$
where the $a_j$'s are $b$-atoms, $\sum_{j}|\lambda_j|<\infty$, and we set
$$
\|f\|_{H^1_{b,{\rm at}}} \colon= \inf\big\{  \sum_{j}|\lambda_j| : \,
\textstyle{ f =\sum_{j}\lambda_j a_j,\, a_j\ b\text{-atoms},\, \sum_{j}|\lambda_j|<\infty } \big\}.
$$
}
\end{definition}
Since  a $b$-atom is in particular a standard atom, if $f
=\sum_{j}\lambda_j a_j \in
H^1_{b,{\rm at}}$,
$\| f\|_{H^1(\bbR)} \le \sum_{j}|\lambda_j|$, which implies that
$\|f\|_{H^1(\bbR)}\le \|f\|_{H^1_{b,{\rm at}}}$.
Moreover,
we have
$\la f,b\ra=0$ for all $f\in H^1_{b,{\rm at}}$, that is, $H^1_{b,{\rm
    at}}\subseteq\ker(b)$. 
Thus, if $f\in
H^1_{b,{\rm at}}$ we have $\|f\|_{H^1_b}\le  \|f\|_{H^1_{b,{\rm at}}}$. 
In \cite{perez1995} 
Perez established the bound 
$$
\sup\left\lbrace \|[b,H](a)\|_{L^1(\bbR)} \colon
  a \text{ is a $b$-atom}\right\rbrace<\infty,
$$
which, in light of the result in \cite{Bownik}, is not sufficient to imply that $[b,H]
\colon H^1_{b,{\rm at}}\to L^1(\bbR)$ is bounded.
L. Ky in \cite{Ky20213} defined another space, that here
we denote by $\cH^1_{b,{\rm at}}$ as
$$
\cH^1_{b,{\rm at}} = \big\{ f\in H^1(\bbR):\, [b,\mathfrak{M}]f \in
  H^1(\bbR)\big\},
  $$
  where $\mathfrak{M}$ is the so-called {\em grand maximal function}.
  Ky 
proved that $[b,T]\colon\cH^1_{b,{\rm at}}\to L^1(\bbR)$ is
bounded, that is $H^1_{b,{\rm at}}\subseteq \cH^1_{b,{\rm at}}$ and that 
$\cH^1_{b,{\rm at}}$ is the largest subspace of $H^1(\bbR)$ on which
$[b,T]$ is bounded.  Again, $T$  denotes a general Cald\'eron--Zygmund operator.

In our opinion, very natural questions to address are the following:\\
For which $b\in \BMO(\bbR)$, we have $H^1_{b, {\rm at}}=H^1_b$?   For
which $b\in \BMO(\bbR)$,  $[b,H]\colon H^1_b\to L^1(\bbR)$?  
}
\end{openproblem}
\medskip

\bibliography{bibliography}
\bibliographystyle{amsplain}

\medskip
\end{document}